\documentclass[12pt]{article}
\usepackage{amssymb}
\usepackage{amssymb}
\usepackage{url}
\usepackage{amsmath, amscd}
\usepackage{latexsym}
\usepackage[dvips]{graphicx}
\usepackage{epstopdf}
\usepackage{lineno}
\usepackage{color}
\usepackage{setspace}
\usepackage{tikz}
\usepackage{tkz-berge}
\usepackage{multicol}
\usepackage{float}
\newtheorem{theorem}{Theorem}[section]
\newtheorem{lemma}[theorem]{Lemma}
\newtheorem{corollary}[theorem]{Corollary}
\newtheorem{definition}[theorem]{Definition}
\newtheorem{proposition}[theorem]{Proposition}

\newtheorem{conjecture}[theorem]{Conjecture}

\newcommand{\qed}{\hfill $\Box$ }
\newcommand{\proof}{\noindent{\bf Proof.}\ \ }
\begin{document}

\title{\Large {\bf Plane Triangulations Without Spanning 2-Trees}}
\author{Allan Bickle\\
\normalsize  Department of Mathematics,\\
\normalsize  Penn State University, Altoona Campus,
\normalsize  Altoona, PA 16601, U.S.A.\\
\small {\tt E-mail: aub742@psu.edu}\\
}

\date{\today}
\maketitle

\begin{abstract}
A 2-tree is a graph that can be formed by starting with a triangle and iterating the operation of making a new vertex adjacent to two adjacent vertices of the existing graph. Leizhen Cai asked in 1995 whether every maximal planar graph contains a spanning 2-tree. We answer this question in the negative by constructing an infinite class of maximal planar graphs that have no spanning 2-tree.

\vskip 0.2in \noindent {\emph{keywords}}:
$k$-tree, triangulation, Hamiltonian cycle, maximal $k$-degenerate graph

\end{abstract}

\section{Introduction}

We consider the problem of whether every maximal planar graph contains
a spanning 2-tree, first proposed by Leizhen Cai \cite{Cai 1995,Cai 1997}
in 1995.
\begin{definition}
A \textbf{$k$-tree} is a graph that can be formed by starting with
$K_{k}$ and iterating the operation of making a new vertex adjacent
to all the vertices of a $k$-clique of the existing graph.
\end{definition}

Note that a $k$-tree is a chordal graph. A more general recursive
construction of $k$-trees is that $K_{k}$ and $K_{k+1}$ are $k$-trees,
and any larger $k$-tree can be formed by identifying two $k$-trees
on $K_{k}$ or $K_{k+1}$.
\begin{definition}
A \textbf{spanning subgraph} of a graph has the same vertex set. A
\textbf{Hamiltonian cycle} of a graph $G$ is a spanning cycle of
$G$. A graph with a Hamiltonian cycle is called a \textbf{Hamiltonian
graph}.

A graph is \textbf{planar} if it has a drawing in the plane that has
no crossings. The \textbf{regions} of a plane drawing are the maximal
pieces of the plane surrounded by edges and vertices. The infinite
region is the \textbf{exterior region}. The \textbf{length} of a region
is the length of a walk around it. A graph is \textbf{maximal planar}
if no edge can be added without making it not planar. A graph is \textbf{outerplanar}
if it has a plane drawing with all vertices on the exterior region.
\end{definition}

Definitions of terms and notation not defined here appear in \cite{Bickle 2020}.
In particular, $K_{n}$ and $C_{n}$ are the complete graph and cycle
of order $n$ and $G+H$ is the join of graphs $G$ and $H$.

Every connected graph has a spanning tree. A spanning 3-tree of a
graph with order $n\geq3$ would have size $3n-6$, so a planar graph
has a spanning 3-tree if and only if it is a 3-tree. The problem of
determining whether there is a spanning 2-tree is more complicated.

Results on $k$-trees and related topics are surveyed in \cite{Bickle 2021}.
Bern \cite{Bern 1987} showed that determining whether a graph has
a spanning $k$-tree is NP-complete when $k\geq2$. Cai and Maffay
\cite{Cai/Maffray 1993} show this is true even for planar graphs
with $\Delta\left(G\right)\leq6$ when $k=2$. Cai found several sufficient
conditions for a spanning 2-tree and showed that it is NP-complete
to determine if $G$ has a spanning $k$-tree even given a spanning
$l$-tree of $G$, $l<k$ \cite{Cai 1995,Cai 1997}. 

Any complete graph has a spanning $k$-tree. Bern \cite{Bern 1987}
showed that it is NP-complete to find a minimum spanning 2-tree for
weighted complete graphs, and found an exponential algorithm for this
problem. Cai \cite{Cai 1995} shows that there is no good approximation
algorithm for weighted complete graphs in general, but there is such
an algorithm when they satisfy the triangle inequality. Shangin and
Pardalos \cite{Shangin/Pardalos 2014} consider various heuristics
for the spanning $k$-tree problem.

Ding \cite{Ding 2016} found applications of spanning $k$-trees to
linguistic grammars, the RNA 3D structure prediction problem, and
learning Markov or Bayesian networks. Spanning 2-trees have applications
in geodesy (geodetic surveying) \cite{McKee 1997} and logic and probability
\cite{Hoover 1990,McKee 2005B}.

\section{Hamiltonian Cycles and 2-Trees}

Leizhen Cai \cite{Cai 1995,Cai 1997} asked in 1995 whether every
plane triangulation contains a spanning 2-tree. Cai did not conjecture
an answer, but I will reframe the problem as a conjecture to simplify
discussion of it.
\begin{conjecture}
\label{conj:Sp2Tree}Every maximal planar graph with order $n\geq3$
contains a spanning 2-tree.
\end{conjecture}

It is easy to show that some special classes of maximal planar graphs
have spanning 2-trees.
\begin{lemma}
\textup{\label{lem:Ham}\cite{Cai 1995} }Every Hamiltonian maximal
planar graph contains a spanning 2-tree.
\end{lemma}

\proof
Adding the edges inside (or outside) a Hamiltonian cycle produces
a spanning 2-tree.
\qed\\

Denote a 4-connected maximal planar graph as a \textbf{4MP}. For these
graphs, the converse is true.
\begin{proposition}
Every spanning 2-tree of a 4MP has a unique Hamiltonian cycle.
\end{proposition}

\proof
A 2-tree has a Hamiltonian cycle if and only if it contains no $K_{2}+\overline{K}_{3}$
\cite{Renjith/Sadagopan 2015}. Let $G$ be a 4MP with a spanning
2-tree $T$. If $T$ is not Hamiltonian, it contains $K_{2}+\overline{K}_{3}$,
so $G$ has a separating triangle and is not 4-connected. Thus $T$
has a Hamiltonian cycle $C$.

A 2-tree with order $n\geq3$ is Hamiltonian if and only if it is
outerplanar \cite{Proskurowski 1979}. It is easily shown by induction
that the cycle is unique, and $T$ can be drawn so that $C$ is the
exterior region.
\qed\\

This shows a correspondence between Hamiltonian cycles and pairs of
spanning 2-trees of 4MPs.
\begin{corollary}
Every 4MP has twice as many spanning 2-trees as Hamiltonian cycles. 
\end{corollary}

Whitney \cite{Whitney 1931} showed that every 4MP is Hamiltonian.
Tutte proved a stronger statement.
\begin{theorem}
\textup{\label{thm:Tutte}\cite{Tutte 1956} Every }planar 4-connected
graph has a Hamiltonian cycle through any two edges of a region.
\end{theorem}

Cai observed the following corollary to Theorem \ref{thm:Tutte}.
\begin{corollary}
\textup{\cite{Cai 1995}} Every 4MP contains a spanning 2-tree.
\end{corollary}

There is another easy special case. Cai and Maffray \cite{Cai/Maffray 1993}
showed that every $l$-tree contains a spanning $k$-tree when $l\geq k\geq1$.
\begin{corollary}
\label{cor:3tree}\textup{\cite{Cai/Maffray 1993} }Every 3-tree contains
a spanning 2-tree.
\end{corollary}

\section{Path-Tree Partitions}

Denote a 4MP or $K_{4}$ as a \textbf{4-block}. Every maximal planar
graph can be formed by identifying triangles of 4-blocks. Any 4-block
has a spanning 2-tree. The question is whether a spanning 2-tree for
the whole graph can be pieced together from those of the 4-blocks.

A spanning 2-tree contains 0, 1, 2, or 3 edges of any given triangle.
If there were some 4MP such that for every spanning 2-tree $T$, there
is some triangle with no edge of $T$, that would be sufficient to
disprove Conjecture \ref{conj:Sp2Tree}. (We could simply attach another
4-block at every triangle, and a spanning 2-tree could not extend
into all of them.) We will show that this is the case.
\begin{definition}
A maximal planar graph has a \textbf{linear Hamiltonian cycle} if
the regions inside (or outside) the cycle share edges with at most
two other such regions (that is, the dual of these regions is a path).
\end{definition}

\begin{conjecture}
\label{conj:LHC}Every 4MP has a linear Hamiltonian cycle.
\end{conjecture}

Conjecture \ref{conj:LHC} is weaker than Conjecture \ref{conj:Sp2Tree}.
We will show that Conjecture \ref{conj:LHC} is false, so Conjecture
\ref{conj:Sp2Tree} is false. To study Conjecture \ref{conj:LHC},
it is convenient to look at the dual graph.
\begin{definition}
The \textbf{Hamiltonian dual} of a planar graph with a given Hamiltonian
cycle is formed by deleting any edges that cross the Hamiltonian cycle
from the dual graph.
\end{definition}

Denote the dual of a 4MP as a \textbf{4MP dual}. A 4MP dual is a 3-connected
cubic planar graph with no nontrivial 3-edge cut. (A trivial edge
cut has all edges incident with a common vertex.)
\begin{definition}
A cubic graph has a \textbf{path-tree partition} if its vertices can
be partitioned into two sets so that one induces a path and the other
induces a tree. A \textbf{path-path partition} and \textbf{tree-tree
partition} are defined similarly. A \textbf{Yutsis graph} is a multigraph
in which the vertex set can be partitioned in two parts such that
each part induces a tree.
\end{definition}

A tree-tree partition is also known as a \textbf{Yutsis decomposition}.
Yutsis graphs have applications in physics, particularly the quantum
theory of angular momenta \cite{YLV 1962}.

A simple degree sum argument shows that the two vertex sets in a tree-tree
partition must have equal size. Theorem \ref{thm:Tutte} implies that
the every 4MP dual has a tree-tree partition.
\begin{proposition}
A 4MP has a linear Hamiltonian cycle if and only if its dual has a
path-tree partition.
\end{proposition}

\proof
$\left(\Rightarrow\right)$ Since all vertices are on a Hamiltonian
cycle, no vertex is inside it. Thus each component of the Hamiltonian
dual is acyclic. Each is clearly connected, so each is a tree. They
have the same order since there are the same number of regions inside
and outside the Hamiltonian cycle. To have a linear Hamiltonian cycle,
one of the trees must be a path.

$\left(\Leftarrow\right)$ If a 4MP dual has a path-tree partition,
the 4MP clearly has a linear Hamiltonian cycle.
\qed\\

To show that Conjecture \ref{conj:LHC} is false, we produce a 4MP
dual that has no path-tree partition.

Let $G_{k}$ have vertices $a_{i}$, $b_{i}$, $c_{i}$, $1\leq i\leq2k$.
The $a_{i}$'s and $b_{i}$'s induce $2k$-cycles, $a_{i}\leftrightarrow c_{i}$,
$b_{i}\leftrightarrow c_{i}$, and $c_{2i-1}\leftrightarrow c_{2i}$
for all $i$, (all mod $2k$). We show $G_{4}$ below.
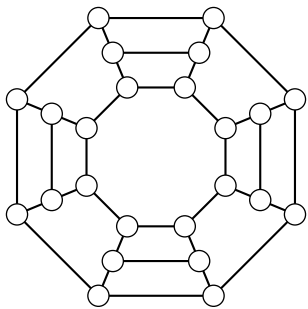
\begin{center}
\begin{tikzpicture}  
\SetVertexSimple[Shape=circle, MinSize=8pt, FillColor=white!50]
\begin{scope}[rotate=22.5]
\grCycle[prefix=a,RA=1]{8}
\grEmptyCycle[prefix=c,RA=1.5]{8}
\grCycle[prefix=b,RA=2]{8}
\end{scope}
\Edges(a0,c0,b0)\Edges(a1,c1,b1)\Edges(a2,c2,b2)\Edges(a3,c3,b3)
\Edges(a4,c4,b4)\Edges(a5,c5,b5)\Edges(a6,c6,b6)\Edges(a7,c7,b7)
\Edges(c1,c2)\Edges(c3,c4)\Edges(c5,c6)\Edges(c7,c0)
\end{tikzpicture}
\par\end{center}
\begin{theorem}
For $k\geq4$, $G_{k}$ has no path-tree partition.
\end{theorem}

\proof
Denote the graph formed from $C_{6}$ by adding a single chord joining
opposite vertices as a \textbf{brick}. Clearly $G_{k}$ contains $k$
bricks. Assume to the contrary that $G_{k}$ has a path-tree partition
with path $P$ and tree $T$. Both $P$ and $T$ must contain at least
one vertex from each cycle.

The ends of $P$ are in one or two bricks, so $P$ must pass through
at least two bricks without ending. It is not possible for $P$ and
$T$ to both pass through the same brick since the two $c$-vertices
would be part of an induced 4-cycle of one of them. If $P$ enters
and exists a brick using two $a$'s, it must end at a $b$ in the
same brick (or vice versa).

To pass through a brick without ending there, $P$ must enter at an
$a$-vertex and exit at a $b$-vertex (or vice versa). Now $P$ must
contain (nonadjacent) $a$-vertices in distinct bricks with a $b$-vertex
between them. But then the graph induced by the vertices not in $P$
is disconnected and hence not a tree, a contradiction.
\qed\\

Essentially the same construction appeared in \cite{BOC 2019}, where
it is used to analyze the number of triangles of certain types in
Hamiltonian maximal planar graphs. Note that $G_{4}$ has order 24.
In fact, a computer search conducted for \cite{BOC 2019} has shown
that the smallest order of a 4MP dual with no path-tree partition
is 24 (personal communication with Gunnar Brinkmann).

The dual of $G_{4}$ is a maximal planar graph with order 14. This
is shown below, where the two black vertices must be identified. Adding
a degree 3 vertex inside each of its 24 regions produces a maximal
planar graph of order 38 with no spanning 2-tree.
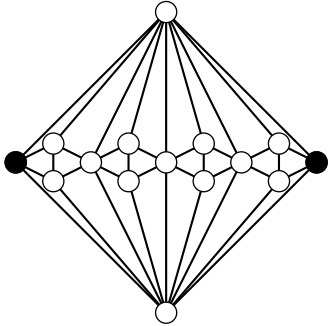
\begin{center}
\begin{tikzpicture}  
\SetVertexSimple[Shape=circle, MinSize=8pt, FillColor=white!50]
\Vertex[x=0,y=1]{01}
\Vertex[x=.5,y=.75]{0a}
\Vertex[x=.5,y=1.25]{0b}
\Vertex[x=1,y=1]{11}
\Vertex[x=1.5,y=.75]{1a}
\Vertex[x=1.5,y=1.25]{1b}
\Vertex[x=2,y=1]{21}
\Vertex[x=2.5,y=.75]{2a}
\Vertex[x=2.5,y=1.25]{2b}
\Vertex[x=3,y=1]{31}
\Vertex[x=3.5,y=.75]{3a}
\Vertex[x=3.5,y=1.25]{3b}
\Vertex[x=4,y=1]{41}
\Vertex[x=2,y=-1]{x}
\Vertex[x=2,y=3]{y}
\AddVertexColor{black}{01,41}
\Edges(01,0a,11,1a,21,2a,31,3a,41,3b,31,2b,21,1b,11,0b,01)
\Edges(x,01,y,0b,0a,x,11,y,1b,1a,x,21,y,2b,2a,x,31,y,3b,3a,x,41,y)
\end{tikzpicture}
\par\end{center}

To produce infinite classes of graphs that do and don't have path-tree
partitions, we need an operation to generate cubic graphs. 
\begin{definition}
Let $uv$ and $wx$ be edges of a cubic graph. Let \textbf{adding
a handle} be the operation of subdividing $uv$ and $wx$ and adding
a new edge $yz$ between the new vertices. Let \textbf{4-handling}
be the operation of adding a handle between two nonadjacent edges
of a region of length 4 of a 4MP dual.
\end{definition}

Every 4MP dual can be constructed from the cube by adding handles
\cite{Faulkner/Younger 1971,Kotzig 1969}.
\begin{proposition}
Let $H$ be formed by 4-handling a 4MP dual $G$. If $G$ has no path-tree
partition then so does $H$.
\end{proposition}

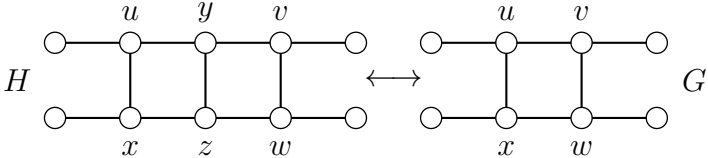
\begin{center}
\begin{tikzpicture}  
\SetVertexSimple[Shape=circle, MinSize=8pt, FillColor=white!50]
\Vertex[x=0,y=0]{00}\node at (-.5,.5) {$H$};
\Vertex[x=0,y=1]{01}
\Vertex[x=1,y=0]{10}\node at (1,-.4) {$x$};
\Vertex[x=1,y=1]{11}\node at (1,1.4) {$u$};
\Vertex[x=2,y=0]{20}\node at (2,-.4) {$z$};
\Vertex[x=2,y=1]{21}\node at (2,1.4) {$y$};
\Vertex[x=3,y=0]{30}\node at (3,-.4) {$w$};
\Vertex[x=3,y=1]{31}\node at (3,1.4) {$v$};
\Vertex[x=4,y=0]{40}
\Vertex[x=4,y=1]{41}
\node at (4.5,.5) {$\longleftrightarrow$};
\Vertex[x=5,y=0]{50}
\Vertex[x=5,y=1]{51}
\Vertex[x=6,y=0]{60}\node at (6,-.4) {$x$};
\Vertex[x=6,y=1]{61}\node at (6,1.4) {$u$};
\Vertex[x=7,y=0]{70}\node at (7,-.4) {$w$};
\Vertex[x=7,y=1]{71}\node at (7,1.4) {$v$};
\Vertex[x=8,y=0]{80}
\Vertex[x=8,y=1]{81}\node at (8.5,.5) {$G$};
\Edges(00,10,20,30,40)\Edges(01,11,21,31,41)
\Edges(10,11)\Edges(20,21)\Edges(30,31)
\Edges(50,60,70,80)\Edges(51,61,71,81)
\Edges(60,61)\Edges(70,71)
\end{tikzpicture}
\par\end{center}
\proof
(contrapositive) Let $uyvwzx$ be a 6-cycle of $H$, and $yz$ be
a handle. Let $uvwx$ be a region of length 4 of $G$ which is 4-handled
to produce the 6-cycle. Assume $H$ has a path-tree partition with
path $P$ and tree $T$. We want to show that $G$ has a path-tree
partition with path $P'$ and tree $T'$.

First suppose that $yz$ is not in $P$ or $T$. Then $y$ and $z$
are not both in $P$ or both in $T$. If $uy$ and $yv$ are both
in $P$ or $T$, let $uv$ be in $P'$ or $T'$, respectively. Similarly,
if $wz$ and $zx$ are both in $P$ or $T$, let $wx$ be in $P'$
or $T'$, respectively. Else don't put $uv$ (or $wx$) in $P'$ or
$T'$. Thus $P'$ and $T'$ are both connected, acyclic, and have
the same order in $G$, so we have a path-tree partition of $G$.

Now suppose that $yz$ is in $T$. At least one of the other vertices,
say $x$, is in $T$. Then $u$ is a leaf of $P$. Then put $ux$
in $T'$ and leave $v$ and $w$ in the same corresponding sets. Thus
$P'$ and $T'$ are both connected, acyclic and have the same order
in $G$, so we have a path-tree partition of $G$. If we exchange
the roles of $P$ and $T$, the argument is similar.
\qed\\

The statement of this proposition does not hold in general for regions
of length more than 4. Since all graphs formed by 4-handling $G_{4}$
have no path-tree partition, we have an infinite class of counterexamples
to Conjecture \ref{conj:Sp2Tree}.

\section{Path-Path Partitions}

The smallest 4MPs are the double wheels $C_{n-2}+\overline{K}_{2}$,
whose duals are the prisms $C_{r}\square K_{2}$, $r\geq4$. Note
that any prism can be generated from the cube by adding handles.
\begin{theorem}
Any 4MP dual constructed from a prism by adding at most two handles
has a path-path partition.
\end{theorem}

Proving this requires many tedious cases. We outline a proof and leave
the details to the reader. A prism $C_{r}\square K_{2}$ has two distinct
regions up to symmetry, a 4-cycle and an $r$-cycle. Denote the edges
joining the two (chordless) $r$-cycles of a prism as \textbf{spokes}.
A prism has a path-path partition using any two spokes and hence omitting
any two given edges of one of the $r$-cycles (see the example below).
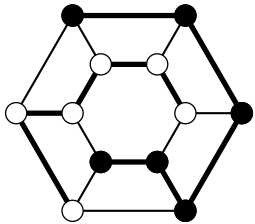
\begin{center}
\begin{tikzpicture}  
\SetVertexSimple[Shape=circle, MinSize=8pt, FillColor=white!50]
\grEmptyCycle[prefix=a,RA=.75]{6}
\grEmptyCycle[prefix=b,RA=1.5]{6}
\AddVertexColor{black}{a4,a5,b5,b0,b1,b2}
\Edges(a5,a0,b0)\Edges(a1,b1)\Edges(a2,b2,b3)
\Edges(a3,a4,b4,b5)
\Edges[local,lw=2pt](a4,a5,b5,b0,b1,b2)
\Edges[local,lw=2pt](a0,a1,a2,a3,b3,b4)
\end{tikzpicture}
\par\end{center}

A handle can always be added to maintain a path-path partition when
both of the edges used are in a path-path partition. One way of adding
a handle to a 4-cycle produces a larger prism and need not be considered
further. The other way produces a path-path partition. The other way
of adding a handle uses two edges of an $r$-cycle, which is always
possible (provided they are nonconsecutive).

When two handles are added, they can be added in two separate 4-cycles
(adjacent or nonadjacent), both in the same 4-cycle (two ways), both
in the same $r$-cycle (independent or not), in two different $r$-cycles
(multiple cases), or one in a 4-cycle and one in an $r$-cycle. In
each case, a path-path partition is easily found.

Note that $G_{4}$ can be generated from a prism by adding four handles.
Next we show that there is a 4MP dual of order 22 with no path-path
partition. It can be formed by adding three handles to a prism. Call
the graph below $H_{22}$. Let the vertices on the exterior 8-cycle
be $a$-vertices, and the vertices on the interior 8-cycle be $b$-vertices.
Call the two edges joining $a$ and $b$-vertices \textbf{spokes}.
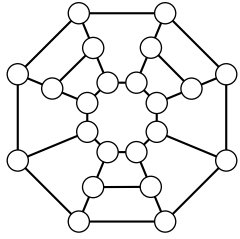
\begin{center}
\begin{tikzpicture}  
\SetVertexSimple[Shape=circle, MinSize=8pt, FillColor=white!50]
\begin{scope}[rotate=22.5]
\grCycle[prefix=a,RA=.5]{8}
\Vertex[x=1,y=0]{c0}
\Vertex[x=.71,y=.71]{c1}
\Vertex[x=0,y=1]{c2}
\Vertex[x=-.71,y=.71]{c3}

\Vertex[x=-.71,y=-.71]{c5}
\Vertex[x=0,y=-1]{c6}

\grCycle[prefix=b,RA=1.5]{8}
\end{scope}
\Edges(a0,c0,b0)\Edges(a1,c1,b1)\Edges(a2,c2,b2)\Edges(a3,c3,b3)
\Edges(a4,b4)\Edges(a5,c5,b5)\Edges(a6,c6,b6)\Edges(a7,b7)
\Edges(c0,c1)\Edges(c2,c3)\Edges(c5,c6)
\end{tikzpicture}
\par\end{center}
\begin{theorem}
The graph $H_{22}$ has no path-path partition.
\end{theorem}

\proof
Assume there is a path-path partition with paths $P_{1}$ and $P_{2}$.
Now $H_{22}$ has three bricks and two spokes. It is not possible
for $P_{1}$ and $P_{2}$ to both pass through the same brick, and
when one does, it must enter at an $a$-vertex and exit at a $b$-vertex
(or vice versa).

If a path passes two bricks, two spokes, or one of each, it either
induces a cycle or disconnects the other path, so this is not possible.
If a path contains a spoke, it could have ends in two bricks, but
would miss the third, so this is not possible.

Thus each path must go through one brick and have ends in the other
two. Thus each path misses all edges between some consecutive pair
of bricks, and this must be a different pair for each path. Thus one
vertex of a spoke cannot be contained in either path, a contradiction.
\qed\\

A computer search by Gunnar Brinkmann shows that 22 is the smallest
order of a 4MP dual with no path-path partition (personal communication).
I had hand-checked that all 4MP duals of order at most 16 have a path-path
partition.

\section{Spanning Maximal 2-degenerate Graphs}

While Conjecture \ref{conj:Sp2Tree} is false, a weaker statement
is true.
\begin{definition}
A graph is \textbf{$k$-degenerate} if its vertices can be successively
deleted so that immediately prior to deletion, each has degree at
most $k$. A graph is \textbf{maximal $k$-degenerate} if no edges
can be added without violating this condition.
\end{definition}

Every $k$-tree is maximal $k$-degenerate, but the converse is false
when $k>1$. We construct a maximal $k$-degenerate graph by starting
with $K_{k}$ and successively adding vertices of degree $k$. Unlike
for a $k$-tree, the neighbors of a new vertex are not required to
induce a clique.
\begin{theorem}
Every maximal planar graph contains a spanning maximal 2-degenerate
graph.
\end{theorem}

\proof
This is obvious for order $n\leq3$. Let $G$ be maximal planar, and
construct it by starting with some 4-block $B_{1}$ and iteratively
adding each new 4-block $B_{r}$ by identifying a triangle $T_{r}$
of $B_{r}$ with a triangle $T_{r}^{*}$ of the existing graph. Let
$G_{1}=B_{1}$ and $G_{r}$ be the graph after $r$ 4-blocks have
been added. We will show that for each $r$, $G_{r}$ has a spanning
maximal 2-degenerate subgraph $M_{r}$.

If $B_{r}$ is a 4-block containing triangle $T_{r}$, then by Theorem
\ref{thm:Tutte} it has a spanning 2-tree that contains $T_{r}$.
Any 2-tree can be constructed starting with any of its triangles.
When identifying $T_{r}$ and $T_{r}^{*}$, delete any edges of $T_{r}$
from the spanning 2-tree of $G_{r}$ that are not in $M_{r-1}$. Thus
the construction of $M_{r-1}$ can continue from $T_{r}$ into $G_{r}$,
producing $M_{r}$. Iterating this process proves the theorem.
\qed\\

\section{Conclusion}

There is more work to do on this problem. In a subsequent paper with
Gunnar Brinkmann, we show that there are maximal planar graphs so
that any 2-tree that is a subgraph of a graph $G$ with order $n$
contains at most $\frac{6n+23}{7}$ vertices. We further examine the
question of what is the smallest order of a maximal planar graph with
no spanning 2-tree. We would like to characterize exactly which maximal
planar graphs have a spanning 2-tree, but that problem seems difficult.

\subsection*{Acknowledgments} 
Thanks to Gunnar Brinkmann for informing me of some relevant articles
and the results of several computer searches of maximal planar graphs
that he conducted.


\begin{thebibliography}{10}
\bibitem{Bern 1987} M. W. Bern, Network design problems: Steiner
trees and spanning k-trees, PhD thesis, University of California,
Berkeley, CA (1987).

\bibitem{Bickle 2020}  A. Bickle, Fundamentals of Graph Theory, AMS
(2020).

\bibitem{Bickle 2021} A. Bickle, A Survey of Maximal k-degenerate
Graphs and k-Trees, 2021+. Submitted.

\bibitem{BOC 2019}  G. Brinkmann, K. Ozeki, and N. Van Cleemput,
Types of triangle in plane Hamiltonian triangulations and applications
to domination and k-walks, Ars Mathematica Contemporanea, 17 1 (2019),
51-66.

\bibitem{Cai 1995} L. Cai, Spanning 2-trees. In: Kanchanasut K.,
Levy JJ. (eds) Algorithms, Concurrency and Knowledge. ACSC 1995. Lecture
Notes in Computer Science, vol 1023. Springer, Berlin, Heidelberg.
(1995) 

\bibitem{Cai 1997} L. Cai, On spanning 2-trees in a graph. Discrete
Appl. Math. 74 (1997), 203\textendash 216.

\bibitem{Cai/Maffray 1993} L. Cai and F. Maffray. On the SPANNING
k-TREE problem. Disc. Appl. Math., 44 (1993), 139-156.

\bibitem{Ding 2016} L. Ding, Maximum Spanning k-Trees: Models and
Applications, PhD Thesis, University of Georgia, 2016.

\bibitem{Faulkner/Younger 1971} G. B. Faulkner, D. H. Younger, The
recursive generation of cyclically k-connected cubic planar maps,
Proceedings of the Twenty-Fifth Summer Meeting of the Canadian Mathematical
Congress, Thunder Bay, 1971, 349\textendash 356.

\bibitem{Hoover 1990} D. R. Hoover, Subset complement addition upper
bounds\textemdash an improved inclusion-exclusion method, J. Statist.
Plann. Inference 24 (1990), 195\textendash 202.

\bibitem{Kotzig 1969} A. Kotzig, Regularly connected trivalent graphs
without non-trivial cuts of cardinality 3, Acta. Fac. Rerum Natur.
Univ. Comenian Math. Publ. 21 (1969), 1\textendash 14.

\bibitem{McKee 1997} T. A. McKee, 2-trees in geodesy, 1850\textendash 1950,
Bull. Inst. Combin. Appl. 21 (1997), 77\textendash 82.

\bibitem{McKee 2005B} T. A. McKee, Spanning (2-)trees of intersection
graphs and Hunter-Worsley-type bounds, Util. Math. 68 (2005), 97-102.

\bibitem{Proskurowski 1979} A. Proskurowski, Minimum dominating cycles
in 2-trees, Internat. J. Comput. Inform. Sci. 8 (1979), 405\textendash 417.

\bibitem{Renjith/Sadagopan 2015} P. Renjith and N. Sadagopan, 2-Trees:
Structural Insights and the study of Hamiltonian Paths, (2015) arXiv:1511.02038

\bibitem{Shangin/Pardalos 2014} R. Shangin and P. Pardalos, Heuristics
for Minimum Spanning k-tree Problem, Procedia Computer Science 31
(2014), 1074-1083.

\bibitem{Tutte 1956} W. T. Tutte, A theorem on planar graphs, Trans.
Amer. Math. Soc. 82 (1956), 99-116.

\bibitem{Whitney 1931} H. Whitney, A theorem on graphs. Ann. of Math.
32 (1931), 378-390.

\bibitem{YLV 1962} A. P. Yutsis, I. B. Levinson, V. V. Vanagas, Mathematical
Apparatus of the Theory of Angular Momentum, Israel Program for Scientific
Translation, Jerusalem, 1962.
\end{thebibliography}
\end{document}